\documentclass[12pt, draftcls, onecolumn]{IEEEtran}
\usepackage{amsmath,amssymb}
\usepackage{graphicx}
\usepackage{verbatim}
\usepackage{subfig}
\usepackage{ifthen}
\usepackage{multirow}
\usepackage{float}
\usepackage{subeqnarray}
\usepackage{enumerate}
\usepackage{cite}
\usepackage[usenames,dvipsnames]{color}

\DeclareMathOperator{\tr}{tr}

\DeclareMathOperator{\rank}{\text{rank }}

\def\ba{\begin{array}}
\def\ea{\end{array}}
\newcommand{\beq}{\begin{equation}}
\newcommand{\eeq}{\end{equation}}
\newcommand{\bq}{\begin{eqnarray}}
\newcommand{\eq}{\end{eqnarray}}
\newcommand{\bqn}{\begin{eqnarray*}}
\newcommand{\eqn}{\end{eqnarray*}}
\newcommand{\bee}{\begin{enumerate}}
\newcommand{\eee}{\end{enumerate}}
\newcommand{\bi}{\begin{itemize}}
\newcommand{\ei}{\end{itemize}}
\newcommand{\ii}{\textbf{i}}

\newcommand{\ve}{\varepsilon}

\renewcommand\Re{\operatorname{\emph{Re}}}

\newtheorem{theorem}{Theorem}
\newtheorem{lemma}[theorem]{Lemma}

\newtheorem{condition}{Condition}

\newtheorem{remark}{Remark}

\newboolean{showcomments}
\setboolean{showcomments}{true}

\newcommand{\bose}[1]{  \ifthenelse{\boolean{showcomments}}
{ \textcolor{Red}{(Bose says:  #1)}} {}  }
\newcommand{\steven}[1]{\ifthenelse{\boolean{showcomments}}
{ \textcolor{Green}{(Steven says: #1)} } {} }
\newcommand{\dennice}[1]{\ifthenelse{\boolean{showcomments}}
{ \textcolor{Green}{(Dennice says:  #1)}}{}}
\newcommand{\mani}[1]{\ifthenelse{\boolean{showcomments}}
{ \textcolor{Red}{(Mani says:  #1)}}{}}

\begin{document}

\title{Quadratically constrained quadratic programs on acyclic graphs
 with application to power flow}

\author{Subhonmesh~Bose,~\IEEEmembership{Student Member,~IEEE,}
Dennice F. Gayme, ~\IEEEmembership{Member,~IEEE,}
K. Mani Chandy,~\IEEEmembership{Fellow,~IEEE,}
and Steven H. Low, ~\IEEEmembership{Fellow,~IEEE}
\thanks{S. Bose is with the Department of Electrical Engineering; S. H. Low and K. M. Chandy are with the Department of Computing and Mathematical Sciences, all at the California Institute of Technology, Pasadena, CA 91125. D. F. Gayme is with the Department of Mechanical Engineering at the Johns Hopkins University, Baltimore, MD 21218
 {\tt \{boses, mani, slow\}@caltech.edu, dennice@jhu.edu}}
\thanks{This work was supported by NSF through NetSE grant CNS 0911041,
Southern California Edison, Cisco, and the Okawa Foundation.}}

\maketitle

\begin{abstract}

This paper proves that non-convex quadratically constrained quadratic programs can be solved in polynomial time when their underlying graph is acyclic, provided the constraints satisfy a certain technical condition. When this condition is not satisfied, we propose a heuristic to obtain a feasible point. We demonstrate this approach on optimal power flow problems over radial networks.

\end{abstract}



\section{Introduction}
\label{sec:intro}

%
%
%

A quadratically constrained quadratic program (QCQP) is an
optimization problem in which the objective function and the constraints
are quadratic. In general, QCQPs are non-convex, and therefore lack computationally efficient solution methods. Many engineering problems
including optimal power flow (OPF) can be
represented as QCQPs with complex variables. The contribution of this paper is to expand the class of non-convex QCQPs for which globally optimal solutions can be guaranteed.

There is a large literature on optimal or approximate algorithms for QCQPs. One such 
method employs a convex semidefinite program that is a rank relaxation of the given QCQP. 
This is commonly referred to as semidefinite relaxation (SDR). Such semidefinite programs 
are solvable in polynomial time using interior-point 
methods \cite{alizadeh1995interior,
  nesterov1987interior, nesterov1988polynomial}. In some 
instances, an optimal solution of the original QCQP can be recovered from an optimal 
solution of its SDR. In other cases, SDR provides a way to approximate the solution of a QCQP. Thus, SDR provides a computationally tractable way to approach QCQPs \cite{luo10, wolkowicz00}. 
For example, SDR has been applied to
a variety of engineering problems such as MIMO antenna beam-forming
\cite{gershman2010convex, sidiropoulos2006transmit, so2008unified,
  huang2010rank}, sensor network localization
\cite{biswas2006semidefinite}, principal component analysis
\cite{mccoy2011two} and stability analysis \cite{liapunov1949probleme}. SDR has also been 
 extensively used in systems and control theory applications \cite{boyd1997linear, boyd1997semidefinite}. Several authors have
investigated exact relaxations, e.g., \cite{kim03, zhang00}, while
others have applied SDR-based approximation techniques to NP-hard
combinatorial problems and non-convex QCQPs, e.g.,
\cite{goemans1995improved, frieze1997improved,
  karger1998approximate}. The accuracy of these approximations has
also been extensively studied, e.g., \cite{ye1999approximating,
  nesterov1997quality, nesterov1998semidefinite}.

In this paper, we prove a sufficient condition under which QCQPs with
underlying acyclic graph structures admit an efficient polynomial time solution using 
its SDR. We then apply our result to the optimal power flow (OPF) problem on radial networks. 
OPF is generally a non-convex, NP-hard problem that 
seeks to minimize some cost function, such as power loss, generation cost and/or 
user utilities, subject to engineering constraints. Since the original formulation of 
Carpentier in 1962 \cite{Carpentier62}, various solution techniques have been used 
for this problem; we refer the reader to \cite{Powerbk, Huneault91,Momoh99a,Momoh99b,Pandya08} 
for some surveys. OPF can be cast as a QCQP. The authors in 
\cite{bai2008, bai2009} propose to solve its SDR as an approximation. The authors in 
\cite{javad10}, instead, propose to solve its (convex) Lagrangian dual 
and provide a sufficient condition under which an optimal solution of the OPF
can be recovered from a dual optimal solution. For IEEE test systems and other randomly generated circuits, 
these approaches have been shown to solve OPF optimally. Recently, OPF 
over radial networks has been of considerable interest. A checkable sufficient condition has been proved in \cite{bose_allerton11, zhang2011geometry, sojoudi2011network} where an optimal solution of OPF can be recovered from an optimal solution of its SDR. This paper extends the previously known class of OPF problems that can be solved efficiently.


The paper is organized as follows. In Section \ref{sec:sdpRelax}, we prove a sufficient condition for
a non-convex QCQP over acyclic graphs to be solvable in polynomial time.
In Section \ref{sec:opf}, we apply this result to characterize a 
class of OPF problems over radial networks that can be solved efficiently. 
In Section \ref{sec:opfNumerics} we describe a heuristic method to obtain 
feasible solutions for QCQPs that do not meet these conditions and thus 
their optimum solution cannot be directly recovered by solving 
its SDR. We further apply this technique for the OPF
problem and demonstrate through simulations that we can always find a near-optimal feasible 
point for OPF. We conclude in Section \ref{sec:conc}.


\section{QCQP's and semidefinite relaxation}
\label{sec:sdpRelax}

Consider the following QCQP with complex variable $x \in
\mathbb{C}^n$. \\
\textbf{Primal problem $P$:}
\bqn
\label{QCQP:P}
\underset{x \in \mathbb{C}^{n}}{\text{minimize}}  & & x^{\mathcal{H}} C  x
\\
\text{subject to:}
& &  x^{\mathcal{H}} C_k x  \leq b_k, \quad k \in \mathcal{K}. \label{eq:P.ineqC}
\eqn
where $x^{\mathcal{H}}$ denotes the conjugate transpose of $x$, $C$ is
either an $n \times n$ complex positive definite matrix (denoted as $C
\succ 0$) or a positive semidefinite matrix (denoted as $C \succeq
0$), $\mathcal{K}$ is a finite index set, $C_{k}$ is an $n \times n$
complex Hermitian matrix and $b_{k}$ is a scalar for each $k \in
\mathcal{K}$.

If the matrices $C_k$, $k \in \mathcal{K}$ are positive semidefinite, then
problem $P$ is a convex program and can be solved in polynomial time \cite{ben2000lectures, boyd04}.
If, however, these matrices are not necessarily positive semidefinite, $P$ is non-convex and NP-hard in general.
The main result of this paper is the identification of a class of QCQPs that can be
solved in polynomial time even though the matrices $C_{k}$, $k \in \mathcal{K}$ are not necessarily
positive semidefinite. This result is applied in the next section
to the optimal power flow problem on radial electric networks.

We begin with some notation.  Let $[n] := \{1, 2, \ldots, n \}$ and
for any matrix $H$, let $H_{ij}$ represent the element corresponding to
the $i^{th}$ row and the $j^{th}$ column.  Define a function
$\mathcal{G}$ from QCQP problems to undirected graphs as follows. For
a QCQP problem $P$, the undirected graph $\mathcal{G}(P)$ has vertex
set $[n]$ and the edge set defined as follows:
\beq
(i, j) \text{ is an edge in } \mathcal{G}(P)
\iff
(i \neq j) \text{ and } (C_{ij} \neq 0 \text{ or } [C_k]_{ij} \neq 0 \text{ for some } k \in \mathcal{K})
\eeq
Since the matrices $C$ and $C_{k}$ are Hermitian, the edges of the
graph are undirected. Also, $\mathcal{G}(P)$ has no self-loops from a
vertex to itself. We restrict attention to QCQP problems $P$ for which
the graph $\mathcal{G}(P)$ is a tree, i.e., it is connected and
acyclic.

\begin{figure*}
        \centering
        \subfloat[Origin does \emph{not} belong to the relative interior of the convex hull of these points.]{
                \includegraphics[width=1.0\textwidth]{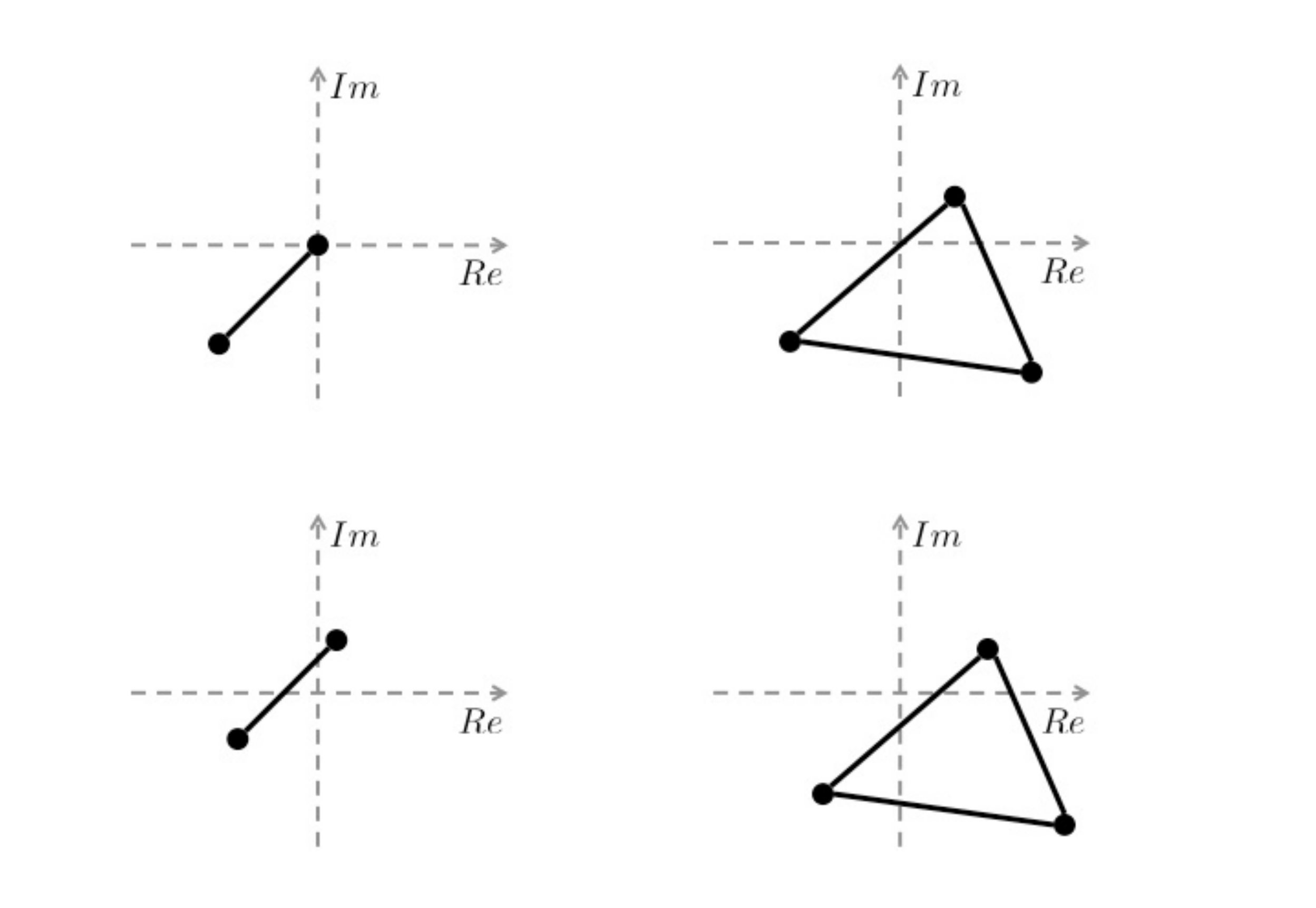}
                \label{fig:allowed}}\\
	\subfloat[Origin lies in the relative interior of the convex hull of these points.]{
                \includegraphics[width=1.0\textwidth]{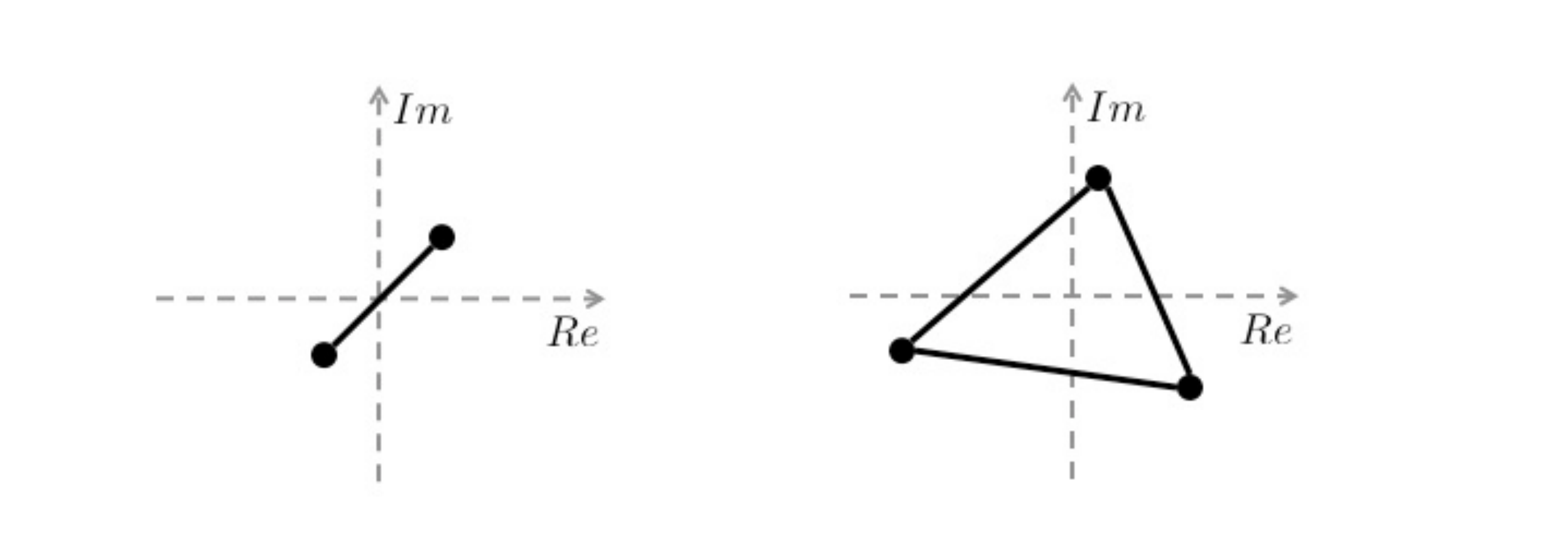}
                \label{fig:notAllowed}}
        \caption{Convex hull of $\{ C_{ij}, [C_k]_{ij}, k \in \mathcal{K} \}$ in condition \ref{ass:qcqpStructure}(\ref{assItem2}).}
\label{fig:hullConditions}
\end{figure*}

For any vector of real numbers $a$, let $a \gg 0$ denote that all
elements of $a$ are strictly positive. For any set of complex numbers
$U = \{ u_1, \ldots, u_r \}$, the \emph{relative interior of the convex hull}
\cite{boyd04} of the set is defined as:
\begin{align}
\label{eq:defRelint}
\{ a_1 u_1 + a_2 u_2 + \ldots a_r u_r \in \mathbb{C} \
\vert \ a \gg 0 \text{ and } \sum_{\ell=1}^{r}a_\ell = 1 \}.
\end{align}
We restrict our attention to QCQP problems for which the origin of the
complex plane does not belong to the relative interior of the convex hull of the set
$\left\{ C_{ij}, [C_{k}]_{i j}, k \in \mathcal{K}\right\}$ for any
edge $(i, j)$ in $\mathcal{G}(P)$.
For any edge $(i, j)$ in $\mathcal{G}(P)$, consider
the points $C_{ij}$ and $[C_{k}]_{i j}, k \in \mathcal{K}$ on the
complex plane. The convex hull of these points is either a single
point, a line segment, or a convex polytope. The condition
states that (a) if the hull is a single point, then that point is not
the origin, (b) if the hull is a line segment, then the origin is
either an extreme point of the line segment or the origin does not lie
on the line segment, and (c) if the hull is a convex polytope then the
origin is either outside or on the boundary of this
polytope. This is illustrated in Figure \ref{fig:hullConditions}. 

We also limit the discussion to QCQPs for which the set of
feasible solutions is bounded and has a strictly
feasible point, i.e., there exists $x \in \mathbb{C}^n$ such that 
$x^\mathcal{H}C_{k} x < b_{k}$ for all $k \in \mathcal{K}$.

To summarize, consider QCQPs that satisfy the
following.
\begin{condition}
\label{ass:qcqpStructure}
\begin{enumerate}[(a)]
\item $\mathcal{G}(P)$ is connected and acyclic.
\label{assItem1}
\item For any edge $(i, j)$ in $\mathcal{G}(P)$, the origin is not in
  the relative interior of the convex hull of $\left\{ C_{ij}, [C_{k}]_{i j}, k
    \in \mathcal{K}\right\}$.
    \label{assItem2}
\item The set of feasible solutions of $P$ is bounded and has a strictly feasible point.
\label{assItem3}
\end{enumerate}
\end{condition}

The main result of the paper is the following theorem and its application.
\begin{theorem}
\label{thm:mainResult}
All QCQPs $P$ that satisfy condition \ref{ass:qcqpStructure} can be
solved in polynomial time.
\end{theorem}
For a continuous optimization problem, we say it can be \emph{solved in polynomial time}
if given any $\zeta > 0$, there is an algorithm that finds a feasible solution to the optimization problem with an
objective value within $\zeta$ of the theoretical optimum in polynomial time \cite{boyd04, ben2000lectures, wolkowicz00}.
For QCQPs $P$ that satisfy condition \ref{ass:qcqpStructure}, Theorem \ref{thm:mainResult} says
that we can construct such a point in polynomial time.

To solve $P$, we use its convex relaxation that can be solved in polynomial time. The relaxation 
is said to be \emph{exact}, if there exists an optimal solution of the
relaxation that can be mapped to an optimal solution of $P$. An exact (convex) relaxation 
by itself does not however guarantee that $P$ can be solved in polynomial time. This is the case 
when the set of optimizers of the relaxation contains solutions that cannot be mapped to a feasible point in $P$ 
and there may or may not be a polynomial time algorithm to find a correct optimum among that set of optimizers.\footnote{
If \emph{every} optimal solution of the relaxation can be mapped to a solution of $P$, then clearly 
$P$ can be solved in polynomial time. This is sufficient but not a necessary condition for a polynomial time solution.}
Theorem \ref{thm:mainResult} asserts that when condition 1 holds, not only the convex relaxation of $P$ is exact, but it 
can also be solved in polynomial time. See Remark 1 at the end of this section for more details.

In the remainder of this section, we prove Theorem \ref{thm:mainResult} for the case where $C$ is
positive definite. The proof for the positive semidefinite case is
presented in the appendix. Our proof requires the following result
on Hermitian matrices, that is presented here and proved in the appendix.
\begin{lemma}
\label{lemma:traceBound}
If $H_1 \succeq 0$ and $H_2 \succeq 0$ are two $n \times n$ matrices,
\bqn
\tr (H_1  H_2) \geq  \rho_{\min}[H_1] \  \rho_{\max}[H_2].
\eqn
where $\rho_{\min}[H]$ and $\rho_{\max}[H]$ respectively denote the minimum and maximum eigenvalues of any Hermitian matrix $H$.
\end{lemma}

\subsection{Proof of Theorem \ref{thm:mainResult} for $C \succ 0$:}
For $C \succ 0$, we prove the result more generally by relaxing the condition that the feasible region of $P$ is bounded. Consider
the following semidefinite program $RP$ where $W$ is an $n\times n$ complex positive
semidefinite matrix.
\\
\textbf{Relaxed Problem $RP$:}
\bq
\label{RP}
\underset{W  \succeq 0}{\text{minimize}}  & & \tr(C W) \nonumber\\
\text{subject to:}
& & \tr \left( C_k W \right)  \leq b_k, \quad k \in \mathcal{K}.
\label{eq:RP.ineqC}
\eq
$RP$ is a convex relaxation of $P$ \cite{wolkowicz00, boyd04}. Define
$p_*$ and $r_*$ as the optimum values of the objective functions for
problems $P$ and $RP$ respectively.
\begin{lemma}
	\label{lemma:relaxedProblem}
	$p_*$, $r_*$ are finite and $p_* \geq r_*$.  If $W_*$ solves $RP$
	optimally and $\text{rank } W_* \leq 1$, then $p_{*} = r_{*}$ and problem $P$ has an exact SDR.
\end{lemma}
\begin{IEEEproof}
	The objective functions of $P$ and $RP$ are nonnegative and
	hence $p_*$ and $r_*$ are finite.
	Given any feasible solution $x$ of $P$, $W := x x^{\mathcal{H}}$ is
	a feasible solution of $RP$. Hence $RP$ is feasible and $p_* \geq
	r_*$.  If $\text{rank } W_{*} = 0$, then
	$W_{*} = 0$, and an optimal solution to $P$ is $x_{*} = 0$, and
	therefore $r_{*} = p_{*}$.  If $\text{rank } W_{*} = 1$ then $W_*$
	has a unique decomposition $W_* = x_*x_*^{\mathcal{H}}$, where
	\bqn
	r_{*} = \tr(C W_{*}) =  x_{*}^{\mathcal{H}} C  x_{*} = p_{*}.
	\eqn
\end{IEEEproof}
Next, we show that there exists a finite $W_*$ that solves $RP$ optimally and has $\rank W_* \leq 1$.



Let the Lagrange multipliers for the inequalities in
\eqref{eq:RP.ineqC} be $\lambda_k \geq 0$ for each $k \in
\mathcal{K}$. Then the Lagrangian dual of $RP$ is
\\
\noindent
\textbf{Dual problem $DP$:}
\bqn
\underset{\lambda \geq 0}{\text{maximize}}
& &   - \displaystyle\sum_{k \in \mathcal{K}} \lambda_k b_k \nonumber
\\
\text{subject to}
& & C + \displaystyle\sum_{k \in \mathcal{K}} {\lambda}_k C_k \succeq 0.
\eqn
For convenience, we introduce the $n \times n$ matrix $A(\lambda)$
defined as:
\beq
A\left({\lambda} \right) := C + \displaystyle\sum_{k \in \mathcal{K}} {\lambda}_k C_k,
\label{eq:defA}
\eeq
Define a function $\mathcal{F}$ from Hermitian matrices to
undirected graphs as follows. For any $n \times n$ Hermitian matrix
$H$, the graph $\mathcal{F}(H)$ on the vertex set $[n]$ satisfies
\beq
(i, j) \text{ is an edge in }
\mathcal{F}(H) \iff i \neq j \text{ and } H_{ij} \neq 0.
\eeq
From the definitions of $\mathcal{F}$, $\mathcal{G}$ and $A(\lambda)$,
it follows that for any $\lambda$, $\mathcal{F}(A(\lambda))$ is a
subgraph of $\mathcal{G}(P)$. For some values of $\lambda$ however,
edge $(i, j)$ may exist in $\mathcal{G}(P)$ but not in
$\mathcal{F}(A(\lambda))$; in this case $\mathcal{F}(A(\lambda))$ is
acyclic but may not be connected, and hence it may be a forest of two
or more disconnected trees rather than a single connected tree that
spans all vertices in the graph. Now, we present a lemma about the
connectedness of the graph $\mathcal{F}(A(\lambda))$.

\begin{lemma}
  	\label{lemma:connectedA}
	For all $\lambda \gg 0$, $\mathcal{F}(A(\lambda))$ is connected.  	
\end{lemma}
\begin{IEEEproof}
	Consider any edge $(i, j)$ in $\mathcal{G}(P)$. From condition
	\ref{ass:qcqpStructure}, the origin is not in the relative interior
	of the convex hull of $\left\{ C_{ij}, [C_{k}]_{ij},
	k \in \mathcal{K}\right\}$. Using \eqref{eq:defA}, we have
	$[A(\lambda)]_{ij} \neq 0$ and hence $(i, j)$ is an edge of
	$\mathcal{F}(A(\lambda))$. Thus $\mathcal{F}(A(\lambda))$ is
	identical to $\mathcal{G}(P)$ which is a tree that spans all the
	vertices of the graph.
\end{IEEEproof}

Next we characterize the relationship between the optimal points of
$RP$ and $DP$. Let $d_*$ denote the optimal value of the objective
function of $DP$.

\begin{lemma}
\label{lemma:rpdpStrongDuality}
	$r_* = d_*$ and $RP/ DP$ has a finite primal dual optimal point $(W_*, \lambda_*)$.
\end{lemma}
\begin{IEEEproof}
	To prove this, we first show that $DP$ is strictly
	feasible. At $\lambda= 0$, $\rho_{\min}[A(\lambda)] = \rho_{\min}[C] > 0$. For a sufficiently small $\lambda' \gg 0$,
	$\rho_{\min}[A(\lambda')] > 0$ and hence $A(\lambda') \succ 0$. This
	implies $\lambda'$ is a strictly feasible point of $DP$. Also since
	$P$ is assumed to be strictly feasible, $RP$ is strictly feasible
	and it has a finite optimum $r_*$. The rest follows from Slater's
	condition \cite{boyd04}.
\end{IEEEproof}
For convenience, define $A_* := A(\lambda_*)$. Now we turn our
attention to the graph of $A_*$, i.e., $\mathcal{F}(A_*)$ to further
analyze the primal dual optimum point $(W_*, \lambda_*)$ of $RP/ DP$.
\begin{lemma}
\label{lemma:tree}
	If $\mathcal{F}(A_*)$ is connected then $\text{rank } W_* \leq 1$.
\end{lemma}
\begin{IEEEproof}
We observe that $\text{rank } A_{*} \geq n - 1$. This follows from a
result in the literature \cite{deVerdiere1998}, \cite[Theorem
3.4]{vanderHolst2003} and \cite[Corollary 3.9]{johnson03} that states that 
for any $n\times n$ positive semidefinite matrix $H$ where 
the associated graph $\mathcal{F}(H)$ is a connected acyclic
graph (i.e., a tree), $\text{rank }H \geq n - 1$.

Next we show that $\text{rank } W_* \leq 1$.  The complementary
slackness condition for optimality of $\left(W_*, {\lambda}_*\right)$
implies
\bqn
\text{tr} (  A_* W_* ) \ = \ 0.
\eqn
Let $W_* = \sum_i \rho_i w_{i} w_{i}^{\mathcal{H}}$ be the spectral
decomposition of $W_*$. Then,
\bqn
\text{tr} (  A_* W_* ) \ = \ \sum_i \ \rho_i \ w_{i}^{\mathcal{H}} A_* w_{i} = 0.
\eqn
Since $A_* \succeq 0$, the eigenvectors $w_{i}$ of $W_*$ corresponding
to nonzero eigenvalues $\rho_i$ are all in the null space of $A_*$.
The rank of $A_*$ is at least $n-1$ and hence its null space has
dimension at most 1, from which it follows that $\text{rank } W_* \leq
1$.
\end{IEEEproof}

$\mathcal{F}(A_*)$ can be connected in one of two ways: (a) the origin of the
complex plane lies strictly outside the convex hull of the set of points
$\{ C_{ij}, [C_k]_{ij} \neq 0, k \in \mathcal{K} \}$ for all edges $(i, j)$ in $\mathcal{G}(P)$, or
(b) $\lambda_* \gg 0$ (from lemma \ref{lemma:connectedA}). In both cases, lemma \ref{lemma:tree} guarantees that $\rank W_* \leq 1$.

If the origin lies on the boundary of the convex hull however, then $\mathcal{F}(A_*)$ may
not be connected when $\lambda_*$ is not element-wise strictly positive and therefore
$\rank W_* \leq 1$ may not hold. In that case, we cannot obtain an optimum solution
of $P$ from the optimum solution of $RP$. We deal with this case where $\mathcal{F}(A_*)$ is
disconnected by using a perturbation \cite{bonnans1998optimization, yildirim2001sensitivity} of $RP/ DP$ so that
$\mathcal{F}(A_*)$ is connected in the perturbed problem.  Then we
recover an optimal solution $W_{*}$ for $RP$ such that $\rank W_* \leq
1$ from the perturbed problem.

Define the perturbed problems for parameter $\ve > 0$:
\\
\noindent
\textbf{Perturbed relaxed problem $RP^{\ve}$:}
\bqn
\underset{W \succeq 0  }
       {\text{minimize}} & & \tr (C W) - \ve \sum_{k \in
    \mathcal{K}} \left[ b_k - \tr \left( C_k W \right) \right]
  \\
  \text{subject to:} & & \tr \left( C_k W \right) \leq {b}_k, \quad k
  \in \mathcal{K}.
\eqn
\noindent
\textbf{Perturbed dual problem $DP^{\ve}$:}
\bqn
\underset{{\lambda}}{\text{maximize}}
& & - \sum_{k \in \mathcal{K}}{\lambda}_k b_k
\nonumber \\
\text{subject to}
& & A\left({\lambda}\right) \succeq 0, \ {\lambda}_{k} \geq \ve, \quad k \in \mathcal{K}.
\eqn
For any variable $z$ in the original problem, let $z^\ve$
denote the corresponding variable in the perturbed problem with
perturbation parameter $\ve$.

\begin{lemma}
\label{lemma:perturbedTree}
There exists a $\ve_0 > 0$, such that for all $\ve$ in $(0, \ve_0)$:
\bee
\item
  $r_*^\ve = d_*^\ve$ and $RP^\ve/ DP^\ve$ has a finite primal
  dual optimal point $(W_*^\ve, \lambda_*^\ve)$.
\item
  $\mathcal{F}(A_*^\ve)$ is connected
  and $\rank W_*^\ve \leq 1$.
\eee
\end{lemma}
\begin{IEEEproof}
  	We choose $\ve_0$ as follows. For $\ve > 0$, observe that
	\bqn
	A(\ve \mathbf{1}) = C + \ve \sum_{k \in \mathcal{K}} C_k,
	\eqn
	where $\mathbf{1}$ is a vector of all ones of appropriate
        size. Then $\rho_{\min}[A(\ve \mathbf{1})]$ is a continuous
        function of $\ve$ that has the value $\rho_{\min}[C] > 0$ at
        $\ve = 0$. Choose $\ve_0$ sufficiently small such that
        $\rho_{\min}[A(\ve \mathbf{1})] $ is strictly bounded away
        from 0, i.e., \bqn \min_{\ve \in [0, \ve_0] }\rho_{\min}[A(\ve
        \mathbf{1})] > 0.  \eqn
	
	Consider any $\ve$ in $(0, \ve_0)$. The feasible sets of $RP$ and $RP^\ve$ are identical. Since $A(\ve \mathbf{1}) \succ 0$ and $W \succeq 0$ for a feasible point $W$ of $RP$ (and $RP^\ve$),
	\begin{align*}
          \tr (C W) - \ve \sum_{k \in \mathcal{K}} \left[ b_k - \tr
            \left( C_k W \right) \right]
		=
           \tr \left[ A(\ve \mathbf{1}) \ W \right] - \ve \sum_{k \in \mathcal{K}} b_k
           	\geq - \ve_{0} \left | \sum_{k \in \mathcal{K}} b_k \right |
	\end{align*}
	and hence $r_*^\ve$ is finite. $RP^\ve/ DP^\ve$ are strictly feasible and $r_*^\ve$ is bounded below.
	The first part of lemma \ref{lemma:perturbedTree} then follows from Slater's condition \cite{boyd04}.

	To prove the second part of lemma \ref{lemma:perturbedTree}, note that
	$\lambda_* \geq \ve \mathbf{1} \gg 0$. Lemmas \ref{lemma:connectedA} and
	\ref{lemma:tree} applied to $RP^\ve$ proves the claim.

\end{IEEEproof}

We have shown that for all $\ve$ in a nonempty interval $(0, \ve_0)$,
the optimal solution $W_*^\ve$ of $RP^\ve$ has rank at most 1. Now we
analyze the behavior of $W_*^\ve$ as $\ve$ converges to zero. Let $\{
\ve_\ell \}_{\ell = 1}^{\infty}$ be a decreasing sequence such that
$\ve_\ell \to 0$ as $\ell \to \infty$. Consider the sequence of
matrices $\{ W_*^{\ve_\ell} \}_{\ell = 1}^{\infty}$; every matrix in
this sequence has rank at most 1. In the next lemma we show that this
sequence has a convergent subsequence and the limit of this
subsequence solves $RP$ optimally.
\begin{lemma}
\label{lemma:convergeW}
$\{ W_*^{\ve_\ell} \}_{\ell = 1}^{\infty}$ has a convergent
subsequence. The limit point $\hat{W}$ of this subsequence solves $RP$
optimally and satisfies $\rank \hat{W} \leq 1$.
\end{lemma}
\begin{IEEEproof}
	Consider any $\ve$ in $(0, \ve_0)$. We first show that $W_*^\ve$ is bounded,
	independent of $\ve$. For any $W$
       in the feasible set of $RP$ (and $RP^\ve$),
	\begin{align*}
	\tr (C W) - \ve \sum_{k \in \mathcal{K}} \underbrace{\left[ b_k - \tr \left( C_k W \right) \right]}_{\geq 0} \leq \tr (C W).
	\end{align*}
	Minimizing both sides over the feasible set of $RP$ (and $RP^\ve$), we have
	\begin{align}
	\label{eq:boundRve}
	r_*^\ve  \leq r_*,
	\end{align}
	that implies
	\begin{align}
	\tr \left[ A(\ve \mathbf{1} ) \ W_*^\ve \right]
	& = r_*^\ve + \ve \sum_{k \in \mathcal{K}} b_k \notag\\
	& \leq r_* + \ve_0 \left\vert \sum_{k \in \mathcal{K}} b_k \right\vert. \label{eq:L8.1}
	\end{align}
	Also, Lemma \ref{lemma:traceBound} implies
	\begin{align}
	\tr \left[ A(\ve \mathbf{1} ) \ W_*^\ve \right]
	& \geq  \ \rho_{\min}[A(\ve \mathbf{1} )] \ \rho_{\max}[W_*^\ve] \notag\\
	& \geq  \ \underbrace{\left[ \min_{\ve' \in [0, \ve_0]}\rho_{\min}[A(\ve' \mathbf{1} )] \right]}_{> 0 \text{ by construction.}} \ \rho_{\max}[W_*^\ve]. \label{eq:L8.2}
	\end{align}
	Combining equations \eqref{eq:L8.1} and \eqref{eq:L8.2}, we obtain a bound for $\rho_{\max}[W_*^\ve]$, independent of $\ve$:
	\begin{align*}
	\rho_{\max}[W_*^\ve]
	\leq
	\frac{r_* + \ve_0 \left\vert \sum_{k \in \mathcal{K}} b_k \right\vert}
	{\min_{\ve' \in [0, \ve_0]}\rho_{\min}[A(\ve' \mathbf{1} )]}.
	\end{align*}

	Thus $W_*^\ve$ is bounded and $\rank W_*^\ve \leq 1$.
	Since the set of positive semidefinite matrices with rank at most $1$ is closed
	\cite{horn2005matrix}, $W_*^\ve$ lies in a compact set and
	$\{ W_*^{\ve_\ell} \}_{\ell = 1}^{\infty}$ has a convergent subsequence.
	Let the limit of this subsequence be $\hat{W}$. Then $\hat{W}$ is feasible for
	$RP$ and $\rank \hat{W} \leq 1$. Next, we prove that $r_* = \tr(C \hat{W})$ and
	hence $\hat{W}$ solves $RP$ optimally.
	
	From \eqref{eq:boundRve},
	\begin{align}
	r_*^\ve = \tr (C W_*^\ve) - \ve \sum_{k \in \mathcal{K}} \left[ b_k - \tr \left( C_k W_*^\ve \right) \right] \leq r_*.
	\end{align}
	Taking limit over the convergent subsequence of $\{ W_*^{\ve_\ell} \}_{\ell = 1}^{\infty}$, we get
	$\tr (C \hat{W}) \leq r_*$. Also, $r_*$ is the optimum value of $RP$ and hence $r_* \leq \tr(C \hat{W})$.
	This completes the proof of lemma \ref{lemma:convergeW}.
\end{IEEEproof}

So far we have shown that $RP$ has a minimizer $\hat{W}$ that satisfies $\rank \hat{W} \leq 1$ and $p_* = r_*$, i.e., SDR of $P$ is exact. But it is, in general, hard to guarantee that solving $RP$ would yield the minimum rank optimizer if the set of optimizers of $RP$ is non-unique. In that case, $RP$ cannot be directly used to solve $P$ in polynomial time. In what follows, we present an algorithm to solve $P$ in polynomial time.

First, solve $RP$ in polynomial time to obtain $r_*$. If the associated optimizer $W_*$ has
rank at most 1, then construct $x_*$ from $W_*$ as in lemma \ref{lemma:relaxedProblem}. We have then found $x_*$ in polynomial time that solves $P$ optimally. If however $\rank W_* > 1$, choose $\ve_0$
as given above and solve $RP^{\ve_0}$ in polynomial time. For any $\ve$ in $(0, \ve_0)$,
\begin{align}
\label{eq:conv.1}
r_*^\ve = \tr (C W_*^\ve) - \ve \sum_{k \in \mathcal{K}}\left[ b_k - \tr (C_k W_*^\ve)  \right] \leq r_*  \leq \tr (CW_*^\ve),
\end{align}
where the first inequality follows from \eqref{eq:boundRve} and the second one follows from the fact that $r_*$ is the optimum value of $RP$. Also, comparing the objective function values of $RP^\ve$ and $RP^{\ve_0}$ at $W_*^\ve$ and $W_*^{\ve_0}$ respectively, we have
\begin{align}
\label{eq:conv.2}
\sum_{k \in \mathcal{K}}\left[ b_k - \tr (C_k W_*^\ve)  \right] \leq \sum_{k \in \mathcal{K}}\left[ b_k - \tr (C_k W_*^{\ve_0})  \right].
\end{align}
Combining \eqref{eq:conv.1} and \eqref{eq:conv.2}, we have
\begin{align*}
\left |  r_* - \tr(CW_*^\ve)  \right |
& \leq \ve \sum_{k \in \mathcal{K}}\left[ b_k - \tr (C_k W_*^{\ve_0})  \right].
\end{align*}
Given any $\zeta > 0$, choose $\ve$ in $(0, \ve_0)$ such that $\sum_{k \in \mathcal{K}}\left[ b_k - \tr (C_k W_*^{\ve_0})  \right] \leq \zeta$. Now solve $RP^\ve$ in polynomial time to get $W_*^\ve$ that satisfies $\rank W_*^\ve \leq 1$ and compute $x_*^\ve$ from it. Then $x_*^\ve$ is a feasible point of $P$ and
\bqn
p_* \leq (x_*^\ve)^\mathcal{H} C (x_*^\ve) \leq p_* + \zeta.
\eqn
Also, we have computed $x_*^\ve$ in polynomial time. This completes the proof of Theorem \ref{thm:mainResult} for the case
where $C$ is positive definite. We remark that a perturbation by an
arbitrary small $\ve>0$ can be represented by treating each perturbed scalar
variable as a pair $[a, a']$ to represent $a + \ve a'$ when solving $RP$ using
any standard polynomial time algorithm like the interior-point method
\cite{alizadeh1995interior, nesterov1987interior, nesterov1988polynomial},
where the pairs $[a, a']$ are ordered lexicographically \cite{dantzig1955generalized}.

The proof extending the theorem to the case where $C$ is positive semidefinite
is given in the appendix.


\begin{remark}
\label{rem:strictFeasible}
	The strict feasibility of $P$ in Condition \ref{ass:qcqpStructure} is required to solve $P$ in polynomial time.
	If we relax that constraint, it can be shown that there still exists a positive semidefinite matrix $W_*$ that
	solves $RP$ optimally and $\rank W_* \leq 1$, i.e., $P$ has an exact SDR. There might also be other optimal solutions of
	$RP$ that do not satisfy the rank condition and hence cannot be mapped to an optimal solution of $P$. 
	Solving for a low-rank optimizer arbitrarily closely 
	in polynomial time is hard to guarantee and is a direction for future work.
\end{remark}


\section{Optimal Power Flow: An application}
\label{sec:opf}

In this section, we apply the results of Section \ref{sec:sdpRelax} to the optimal power flow (OPF) problem.
We start by summarizing some of the recent results on OPF relaxations in Section  \ref{subsec:pw}. In
Section \ref{subsec:pf} we formulate OPF as a QCQP. In Section \ref{subsec:opfResults} we restrict our attention to OPF over radial networks, which are the networks commonly found in distribution circuits, and use Theorem \ref{thm:mainResult} to characterize a set of conditions under which OPF can be solved efficiently.

\subsection{Prior work}
\label{subsec:pw}
As previously discussed, OPF can be cast as a QCQP. Various non-linear programming techniques have been applied to the resulting nonconvex problem, e.g., in \cite{torres1998interior, nejdawi2000efficient, jabr2003primal}. An SDR for OPF has been explored in \cite{bai2008, bai2009} and their simulations indicate that the SDR provides an exact solution of original OPF for many of the IEEE test systems \cite{UW_data}. The authors in \cite{javad10, Lavaei2011} propose to solve the convex Lagrangian dual of OPF and derive a sufficient condition under which an optimal solution of OPF can be recovered from an optimal
dual solution.
 Though SDR recovers a solution to OPF on IEEE test systems, it does not work on all problem instances.
 Such limitations have been most recently reported
in \cite{Lesieutre-2011-OPFSDP-Allerton}, though the nonconvexity of power flow solutions have been
studied much earlier, e.g., in
\cite{Hiskens-2001-OPFboundary-TPS, LesieutreHiskens-2005-OPF-TPS, Hill-2008-OPFboundary-TPS, zhang2011geometry}. 

Recently, a series of work has explored a class of problems where such limitations do not apply due to the network topology. It has been independently reported in \cite{bose_allerton11, zhang2011geometry, sojoudi2011network} that the SDR of OPF is
exact for radial networks provided certain conditions on the power flow constraints are satisfied. A different approach to OPF has been explored using the branch flow model, first introduced in \cite{Baran1989a, Baran1989b}. While \cite{Taylor2011PhD} studies a linear approximation of this model, various relaxations based on second-order cone programming (SOCP) have been proposed in \cite{Jabr2006, Masoud2011, Gan-2012-BFMt, Na-2012-BFMt}. Authors in \cite{Masoud2011, Gan-2012-BFMt, Na-2012-BFMt} prove that this relaxation is exact for radial networks when there are no upper bounds on loads, or when there are no upper bounds on voltage magnitudes. 

Motivated by the results in \cite{Masoud2011}, a more general branch flow model is introduced in
\cite{Farivar-2012-BFM-CDC} for the power flow analysis and optimization for both radial and meshed
networks.     The precise relationships between the SOCP relaxations and the SDR for the OPF problem has been recently identified in \cite{bose2012equivalence}.


\subsection{Problem Formulation}
\label{subsec:pf}

Consider a power system network with $n$ nodes (buses).
The admittance-to-ground at bus $i$ is $y_{ii}$ and the
admittance of the line between connected nodes $i$ and $j$ (denoted by $i\sim j$) is $y_{ij} = g_{ij} - \ii b_{ij}$.
We assume both $g_{ij} > 0$ and $b_{ij} > 0$, i.e., the lines
are resistive and inductive.  Define the corresponding $n \times n$ admittance matrix $Y$ as
\begin{align}\label{eq:defY}
Y_{ij} = \begin{cases} y_{ii} + \displaystyle\sum_{j \sim i} y_{ij}, & \text{ if } i=j,\\
- y_{ij}, & \text{ if } i \neq j \text { and $i \sim j$}, \\
0 & \text{ otherwise}.\end{cases}
\end{align}
\noindent\begin{remark}
$Y$ is symmetric but not necessarily Hermitian.
\end{remark}

The remaining circuit parameters and their relations are defined as follows.
\begin{itemize}
\item $V$ and $I$ are $n$-dimensional complex voltage
and current vectors, where $V_k$, $I_k$ denote the voltage and the injection current
at bus $k \in [n]$ respectively.  The voltage magnitude at each bus is bounded as
\bqn
0 < \underline{W}_k \leq |V_k|^2 \leq \overline{W}_k, \,\quad k\in[n].
\eqn
\item $S = P + \ii Q$ is the $n$-dimensional complex power vector, where
	 $P$ and $Q$ respectively denote the real and reactive powers and	
\begin{equation}
\label{eqn:Sk}
S_k  = P_k + \ii Q_k =  V_k I_k^{\mathcal{H}}, \;k \in[n].
\end{equation}
\item $P_k^{D}$ and $Q_k^{D}$ are the real and reactive power demands
at bus $k\in[n]$, which are assumed to be fixed and given.

\item $P_k^{G}$ and $Q_k^{G}$ are the real and reactive power generation
at bus $k$.   They are decision variables that satisfy the constraints
$\underline{P}_k^G \leq P_k^G \leq \overline{P}_k^G$ and
$\underline{Q}_k^G \leq Q_k^G \leq \overline{Q}_k^G$.

\end{itemize}
Power balance at each bus $k\in[n]$ requires $P_k^G = P_k^D + P_k$ and
$Q_k^G = Q_k^D + Q_k$, which leads us to define
\bqn
\underline{P}_k \ := \ \underline{P}_k^G -P_k^{D}, & &
\overline{P}_k  \ := \ \overline{P}_k^{G} -P_k^{D}
\\
\underline{Q}_k \ := \ \underline{Q}_k^G -Q_k^{D}, & &
\overline{Q}_k  \ := \ \overline{Q}_k^{G} -Q_k^{D}.
\eqn
The power injections at each bus $k\in[n]$ are then bounded as
\bqn
\underline{P}_k  \leq  P_k  \leq  \overline{P}_k,
& &
\underline{Q}_k  \leq  Q_k  \leq  \overline{Q}_k.
\eqn
The branch power flows and their limits are defined as follows.
\begin{itemize}
\item $S_{ij} = P_{ij} + \ii Q_{ij}$ is the sending-end complex power flow from node $i$ to node $j$,
where $P_{ij}$ and $Q_{ij}$ are the real and reactive power flows respectively.  The real power flows
are constrained as $| P_{ij} | \leq \overline{F}_{ij}$ where $\overline{F}_{ij}$ is the line-flow limit
between nodes $i$ and $j$ and $\overline{F}_{ij} = \overline{F}_{ji}$.

\item $L_{ij} = P_{ij} + P_{ji}$ is the power loss over the line between nodes $i$ and $j$, satisfying
$L_{ij} \leq \overline{L}_{ij}$ where $\overline{L}_{ij}$ is the thermal line limit and $\overline{L}_{ij} = \overline{L}_{ji}$. Also, observe that since $L_{ij} \geq 0$, we have $|P_{ij} | \leq \overline{F}_{ij},  |P_{ji} | \leq \overline{F}_{ji}$ if and only if $P_{ij} \leq \overline{F}_{ij},   P_{ji} \leq \overline{F}_{ji} $.

\end{itemize}

Let $J_k = e_k e_k^{\mathcal{H}}$ where $e_k$ is the $k$-th canonical basis vector in $\mathbb{C}^n$. Define $Y_k := e_k e_k^{\mathcal{H}} Y$. Substituting these expressions into \eqref{eqn:Sk} yields
\bq
S_k
&= & e_k^{\mathcal{H}} V I^\mathcal{H} e_k = \tr \left( VV^\mathcal{H} (Y^\mathcal{H} e_k e_k^{\mathcal{H}}) \right) =  V^\mathcal{H} Y_k^\mathcal{H} V
\nonumber\\
&= & \left( V^\mathcal{H}  \underbrace{\left(\frac{Y_k^\mathcal{H} + Y_k}{2}\right)}_{=:\Phi_k} V \right) +
\ii  \left( V^\mathcal{H} \underbrace{\left(\frac{Y_k^\mathcal{H} - Y_k}{2\ii}\right)}_{=:\Psi_k} V \right),
\label{eq:defPhiPsi}
\eq
where $\Phi_k$ and $\Psi_k$ are Hermitian matrices.  Thus,
the two quantities $V^\mathcal{H} \Phi_k V$ and $V^\mathcal{H} \Psi_k V$ are real numbers;
moreover
\begin{align*}
P_k
=  V^\mathcal{H} \Phi_k V, \quad
Q_k
=  V^\mathcal{H} \Psi_k V.
\end{align*}
The real power flow from $i$ to $j$  can be expressed as a quadratic form as follows.
\beq
P_{ij}    =  \Re \{ V_i (V_i - V_j)^{\mathcal{H}} y_{ij}^{\mathcal{H}} \} = V^{\mathcal{H}} M^{ij} V, \label{eq:defPij}
\eeq
where $M^{ij}$ is an $n \times n$ Hermitian matrix. Further details of the OPF problem formulation are provided in the appendix.

The thermal loss of the line connecting buses $i$ and $j$ is
\beq
\label{eq:defTij}
L_{ij} = L_{ji} = P_{ij} + P_{ji} = V^{\mathcal{H}} T^{ij} V
\eeq
where $T^{ij} = T^{ji} := M^{ij} + M^{ji} \succeq 0$.

For a Hermitian positive semidefinite $n \times n$ matrix $C$, we have \\
\noindent
\textbf{Optimal power flow problem $OPF$:}
\begin{subequations}
\label{OPF}
\begin{align*}
& \underset{V \in \mathbb{C}^n}{\text{minimize}}   \; \; V^\mathcal{H} C V &\\
& \text{subject to:} &
\end{align*}
\begin{align}
& \underline{P}_k \leq  V^\mathcal{H} \Phi_k V  \leq \overline{P}_k, \quad k \in [n],
\label{eq:opf.P}\\
& \underline{Q}_k \leq V^\mathcal{H} \Psi_k V  \leq \overline{Q}_k, \quad k \in [n],
\label{eq:opf.Q}\\
& \underline{W}_k \leq V^\mathcal{H} J_k V   \leq \overline{W}_k, \quad k \in [n],
\label{eq:opf.V}\\
& V^\mathcal{H} M^{ij} V   \leq \overline{F}_{ij}, \quad i \sim j,
\label{eq:opf.F}\\
& V^\mathcal{H} T^{ij} V   \leq \overline{L}_{ij}, \quad i \sim j,
\label{eq:opf.L}
\end{align}
\end{subequations}
where \eqref{eq:opf.P}--\eqref{eq:opf.L} are respectively constraints on the real and reactive powers, the voltage magnitudes, the
line flows and thermal losses.
Note that since $T^{ij} \succeq 0$, \eqref{eq:defTij} implies that $P_{ij}+P_{ji} \geq 0$.
This means that \eqref{eq:opf.F} holds if and only if $|P_{ij}| \leq \overline{F}_{ij}$, i.e., \eqref{eq:opf.F}
bounds the line flows on both ends.

We do not include line-flow constraints that impose an upper bound on the apparent power $\sqrt{P_{ij}^2 + Q_{ij}^2}$
on each branch  $i \sim j$. These constraints are not quadratic in voltages and hence beyond the scope of our model.

\begin{remark}[Objective Functions]
\label{rem:objFunction}
We consider different optimality criteria by changing $C$ as follows:
\begin{itemize}
\item Voltages: $C={\mathcal{I}}_{n \times n}$ (identity matrix) where we aim to minimize $\| V \|^2 = \sum_k | V_k |^2$.
\item Power loss: $C = (Y + Y^{\mathcal{H}})/2$ where we aim to minimize
	$\sum_i g_{ii}|V_i|^2 + \sum_{i<j} g_{ij} |V_i - V_j|^2$.
\item Production costs: $C = \sum_{k} c_k \Phi_k$ where we aim to minimize $\sum_k c_k P_k^G$, $c_k \geq 0$.
\end{itemize}
We assume $C$ is positive semidefinite.
\end{remark}

\subsection{Semidefinite relaxation of OPF over radial networks}
\label{subsec:opfResults}

Assume hereafter that $OPF$ is feasible.  To conform to the notations of Section \ref{sec:sdpRelax}, replace
the constraint in \eqref{eq:opf.P} by the equivalent constraints
\bqn
V^\mathcal{H}[ \Phi_k ] V & \leq & \overline{P}_k, \quad k \in [n] \\
V^\mathcal{H}[- \Phi_k] V & \leq & -\underline{P}_k, \quad k\in [n].
\eqn
Similarly rewrite \eqref{eq:opf.Q} and \eqref{eq:opf.V}. Then the set of matrices $\{C_k, k\in \mathcal{K}\}$ and the set of scalars $\{b_k, k \in \mathcal{K} \}$ in $OPF$ are defined as
\begin{align}
\left\{ C_k, k\in \mathcal{K} \right\}
& := \left\{\Phi_k, -\Phi_k, \Psi_k, -\Psi_k, J_k, -J_k, \ k \in [n]\right\} \ \ \bigcup \ \ \left\{M^{ij}, T^{ij}, \  i \sim j\right\}\\
\{b_k, k \in \mathcal{K} \}
& := \left\{\overline{P}_k, -\underline{P}_k, \overline{Q}_k, -\underline{Q}_k, \overline{W}_k, -\underline{W}_k, \ k \in [n]\right\} \ \ \bigcup \ \ \left\{\overline{F}^{ij}, \overline{L}^{ij}, \  i \sim j\right\}
\end{align}

We limit the discussion to $OPF$ instances where the graph of the power network is a tree $\mathcal{T}$ on $n$ nodes. It can be checked that the graph of the problem $OPF$ satisfies
\beq
\mathcal{G} (OPF) = \mathcal{T}.
\eeq
Thus condition \ref{ass:qcqpStructure}(\ref{assItem1}) holds for $OPF$ over $\mathcal{T}$.
In general, condition \ref{ass:qcqpStructure}(\ref{assItem2}) does not hold for $OPF$ over $\mathcal{T}$. To illustrate this point, consider an edge $(i, j)$ in $\mathcal{T}$. The admittance of the line joining buses $i$ and $j$ is $g_{ij} - \ii b_{ij}$. Then $[C_k]_{ij}, k \in \mathcal{K} $ are given as (details are in the appendix):
\begin{enumerate}[(a)]
\item $[\Phi_i]_{ij} = -g_{ij}/2 + \ii b_{ij}/2$,
\item $[\Phi_j]_{ij} = -g_{ij}/2 - \ii b_{ij}/2$,
\item $[\Psi_i]_{ij} = -b_{ij}/2 - \ii g_{ij}/2$,
\item $[\Psi_j]_{ij} = -b_{ij}/2 + \ii g_{ij}/2$,
\item $[M^{ij}]_{ij} =  -g_{ij}/2 + \ii b_{ij}/2$,
\item $[M^{ji}]_{ij} =  -g_{ij}/2 - \ii b_{ij}/2$,
\item $[T^{ij}]_{ij} = [T^{ji}]_{ij} =  - g_{ij}$.
\end{enumerate}
For the objective functions considered, we have
\begin{enumerate}[(a)]
\item Voltages: $C_{ij} = 0$,
\item Power loss: $C_{ij}  = -g_{ij}$,
\item Production costs: $C_{ij} =  -g_{ij} (c_i + c_j)/2 + \ii b_{ij} (c_i - c_j)/2$.
\end{enumerate}
\begin{figure}
	\centering
  	\includegraphics[width=0.8\textwidth]{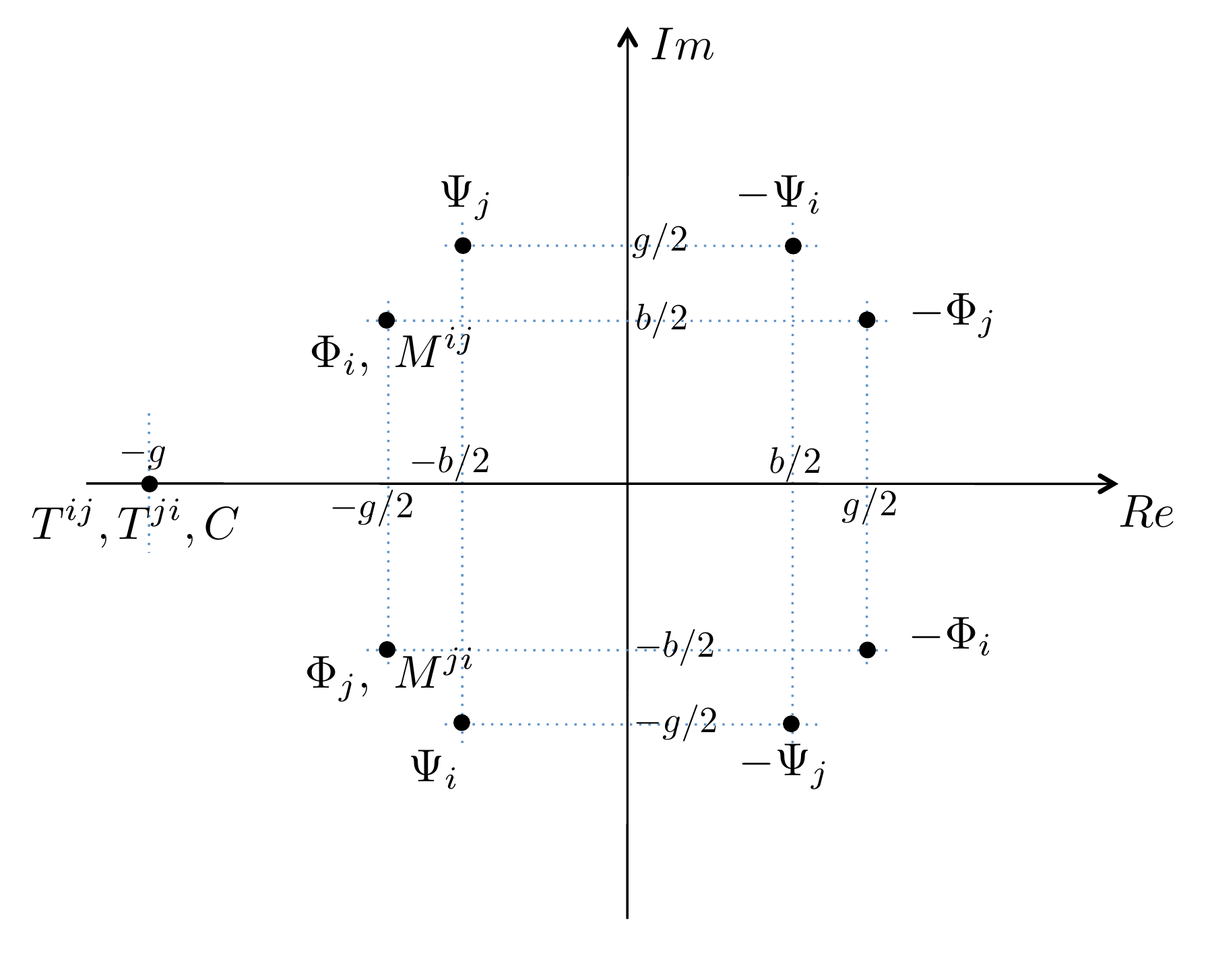}
	\caption{$C_{ij}$ and non-zero $\left( [C_k]_{ij}, k \in \mathcal{K} \right)$ on the complex plane for $OPF$ for a fixed line $(i, j)$ in tree $\mathcal{T}$.}
  	\label{fig:opfConditions}
\end{figure}
For the purpose of this discussion, consider power-loss minimization as the objective, i.e., $C_{ij} = -g_{ij}$. Also, assume $g_{ij} > b_{ij} > 0$. We plot the non-zero $(i,j)$-th entries of the matrices  $C_k, k \in \mathcal{K}$ and $C$,  on the complex plane in Figure \ref{fig:opfConditions} and label each point with its corresponding matrix. Clearly if we consider all the points $[C_k]_{ij}, k \in \mathcal{K}$ and $C_{ij}$, the origin of the complex plane lies in the relative interior of the convex hull of these points, i.e., condition \ref{ass:qcqpStructure}(\ref{assItem2}) does not hold.

To apply Theorem \ref{thm:mainResult} to $OPF$, consider the index-set $\tilde{\mathcal{K}} \subseteq \mathcal{K}$ so that condition \ref{ass:qcqpStructure}(\ref{assItem2}) holds for $C_k, k \in \tilde{\mathcal{K}}$ and $C$. This corresponds to removing certain inequalities in $OPF$, i.e., $b_k = + \infty$ for $k \in \mathcal{K} \setminus \tilde{\mathcal{K}}$. For example, removing $-\Phi_j$ from the set $\{C_k, k \in \mathcal{K}\}$ corresponds to setting $\underline{P}_j = -\infty$.

Condition \ref{ass:qcqpStructure}(\ref{assItem3}) requires that the feasible set of $OPF$ be bounded.
This is always the case when $\overline{W}_k$ is finite for all buses $k \in [n]$.
We discuss the technical condition of strict feasibility in remark \ref{rem:strictFeasibleOPF}.

\begin{theorem}
\label{thm:opfResult}
For $\tilde{\mathcal{K}} \subseteq \mathcal{K}$, suppose condition \ref{ass:qcqpStructure} holds for $OPF$ with $C_k, k \in \tilde{\mathcal{K}}$ and $C$. Then $OPF$ can be solved in polynomial time.
\end{theorem}
\begin{remark}
\label{rem:strictFeasibleOPF} When $OPF$ is strictly feasible, it can be
solved in polynomial time. When strict feasibility does not hold, it follows
from remark \ref{rem:strictFeasible} that $OPF$ still has an exact SDR but we do not guarantee a polynomial time solution.
\end{remark}

We explore, through examples, some constraint patterns for $OPF$ over radial networks where a polynomial time solution (for strictly feasible OPF instances) or an exact SDR (for OPFs that may not be strictly feasible) is guaranteed.\\

\emph{Example 1:} In Figure \ref{fig:opfConditions}, consider the $(i, j)$-th elements of the following set of matrices:
\bqn
	\left\{ \Phi_i,  \Phi_j,  \Psi_i, \Psi_j, -\Psi_i, M^{ij}, M^{ji}, T^{ij} = T^{ji}, C \right\}.
\eqn
The origin of the complex plane lies on the boundary of the convex hull of these points but not in its relative interior. With this set of points, associate a constraint pattern defined as follows. For any point in the diagram that is not a part of this set, the inequality associated with that matrix is removed from $OPF$. For example, the matrices $-\Phi_j$, $-\Phi_i$ and $-\Psi_j$ do not feature on the list of points. Hence,
\beq
\label{eq:constPattern1}
\underline{P}_j = \underline{P}_i = \underline{Q}_j = -\infty.
\eeq
This can be generalized to a constraint pattern over $\mathcal{T}$ by removing the lower bounds on real powers at all nodes and the lower bounds on reactive powers at alternate nodes.

\emph{Example 2:}
Suppose $\underline{P}_{k} = \underline{Q}_{k} = -\infty$ for all nodes $k$ in $\mathcal{T}$. This corresponds to considering points only on the left-half plane in Figure \ref{fig:opfConditions} for all edges $(i, j)$ in $\mathcal{T}$. Clearly, condition \ref{ass:qcqpStructure}(\ref{assItem2}) holds in this case. In Figure \ref{fig:opfConditions}, we assume $g_{ij} > b_{ij} > 0$. Regardless of the ordering between $g_{ij}$ and $b_{ij}$ for edges $(i, j)$ in $\mathcal{T}$, the set of points considered in this constraint pattern always lies in the left half of the complex plane.

Removing the lower bounds in real and reactive power can be interpreted  as load over-satisfaction, i.e., the real and reactive powers supplied to a node $k$ can be greater than their respective real and reactive power demands $P_k^D$ and $Q_k^D$. $OPF$ on a radial network with load over-satisfaction can be solved efficiently. This result has been reported in \cite{bose_allerton11, zhang2011geometry, sojoudi2011network}.

\emph{Example 3:} Consider voltage minimization, i.e., $C = \mathcal{I}_{n \times n}$. In Figure \ref{fig:opfConditions}, consider the $(i,j)$-th entries of the following set of matrices:
\bqn
	\left\{ -\Phi_i, \Phi_j,  -\Phi_j,  \Psi_i, -\Psi_j, C \right\}.
\eqn
The origin of the complex plane is again on the boundary of the convex hull of these points. The constraint pattern associated with this set of points is
\bqn
\overline{P}_i =  \overline{Q}_j = \overline{L}_{ij} = \overline{L}_{ji} = \overline{F}_{ij} = \overline{F}_{ji} = + \infty,
\quad \mbox{and} \quad
\underline{Q}_i = -\infty.
\eqn
This constraint pattern is consistent with condition \ref{ass:qcqpStructure}(\ref{assItem2}) over the edge $(i, j)$ and we can construct a constraint pattern for the OPF problem.

\section{Numerical examples}
\label{sec:opfNumerics}


\subsection{Numerical techniques}
\label{subsec:sdpNumerics}

In Section \ref{sec:sdpRelax}, we have identified conditions under which an SDR of a QCQP over an acyclic connected graph can be used to solve the nonconvex problem $P$ efficiently. When these conditions are not satisfied, the QCQP may not be polynomial time computable or have an exact SDR, i.e. $RP$ yields an optimal $W_*$ with $\rank W_*>1$. In that case, $W_*$ cannot be mapped to an optimal $x_*$ for the problem $P$. In this section we propose a method to construct a feasible solution $\tilde{x}$ for $P$ using such an optimal $W_*$ of $RP$. The following relation characterizes the relationship between the optimal solution of $P$ and its value at $\tilde{x}$:
\bq
\mbox{ objective value of $RP$ at $W_*$} & \leq & \mbox{optimum objective value of $P$}\nonumber\\
& \leq & \mbox{objective value of $P$ at $\tilde{x}$}.
\label{eq:justNR}
\eq
In many practical problems where $\text{rank } W_{*} > 1$, the principal eigenvalue of $W_*$ is
orders of magnitude greater than the other eigenvalues. We use the principal eigenvector to search for a ``nearby'' feasible point of $P$ as follows. Let $w_* \in \mathbb{C}^n$ be the principal eigenvector of $W_*$ and define the starting point $x_0$ of the algorithm as
\bqn
 x_0 := w_* \sqrt{\tr (C W_*)}.
\eqn
This scaling ensures that the objective value at $x_0$ is the optimum objective value of $RP$ at $W_*$. If $x_0$ satisfies all constraints in $P$, then the algorithm ends with $\tilde{x} = x_0$. Otherwise, we construct a sequence of points $(x_1, x_2, \ldots)$ where $x_{m+1}$ is constructed from $x_m$ as follows.
\begin{enumerate}
\item
For $k \in \mathcal{K}$, linearize the function $f_k (x) = x^{\mathcal{H}} C_k x$ around the point $x_m$ and call this function $f_k^{(m)} (x)$, i.e.,
\bqn
f_k^{(m)} (x) = x_m^{\mathcal{H}} C_k x_m + 2 \Re \left[ x_m^{\mathcal{H}} C_k (x - x_m ) \right].
\eqn
\item
For $k \in \mathcal{K}$, define
\bqn
s_k^{(m)} (x) :=
\begin{cases}
\underline{b}_k - f_k^{(m)}(x), & \text{if } f_k^{(m)}(x) \leq \underline{b}_k,\\
0 & \text{if } \underline{b}_k \leq f_k^{(m)}(x) \leq \overline{b}_k,\\
 f_k^{(m)}(x) - \overline{b}_k, & \text{if } \overline{b}_k \leq f_k^{(m)}(x).
 \end{cases}
\eqn
We can interpret $s_k^{(m)}(x)$ as the amount by which the linearized function $f_k^{(m)}$ violates the inequality constraint $\underline{b}_k \leq f_k^{(m)}(x) \leq \overline{b}_k$.

\item
Compute $x_{m+1}$ using
\bqn
\label{iterationProb}
x_{m+1} = \underset{x \in \mathbb{C}^{n}}{\arg \min}  & & \sum_{k \in \mathcal{K}} [s_k^{(m)} (x)]^2\\
\text{subject to:}
& & \| x - x_m \|_1 \leq \gamma, \label{eq:stepSize}
\eqn
where $\| . \|_1$ denotes the $\ell^1$ norm and $\gamma$ is the maximum allowable step-size. This is a parameter for the algorithm and should be chosen such that the linearization $f_k^{m} (x) $ is a reasonably good approximation of the quadratic form $f_k(x)$ for all $k \in \mathcal{K}$ in the $\ell^{1}$ ball centered around $x_m$ with radius $\gamma$.
\item
If $x_{m+1}$ satisfies all constraints in $P$, then the algorithm ends with $\tilde{x} = x_{m+1}$.
\end{enumerate}
This heuristic approach either ends at a feasible point $\tilde{x}$ of $P$ or it fails to produce one within a fixed number of iterations. In the next section, we show that this technique performs quite well for numerical $OPF$ examples where the SDR yields an optimal solution $W_*$ with rank more than 1.


\subsection{OPF test examples}

The SDR of $OPF$ and the techniques described in section \ref{sec:opf} are illustrated on a sample distribution circuit from Southern California and randomly generated radial circuits. The semidefinite program is solved in MATLAB using YALMIP \cite{lofberg2004yalmip}. If the solution yielded $W_*$ such that $\text{rank } W_* \leq 1$ then the optimal voltage profile ($V_*$) to the $OPF$ problem is calculated from $W_* = (V_*)(V_*)^{\mathcal{H}}$. If the optimal $W_*$ does not satisfy the rank condition, the heuristic approach described above is used to find a feasible point of $OPF$. The feasible point obtained may not be optimal for the original problem, so we characterized its sub-optimality by defining the following quantity.
\bqn
\eta := \frac{\text{Objective value at heuristically reached feasible point}}{\text{Objective value at optimal point of relaxed problem}} - 1.
\eqn
Smaller values of $\eta$ indicate that the feasible point obtained using the algorithm is close to the theoretical optimum of $OPF$.

Throughout this section, let $y = (a, b)$ denote a $y$ drawn from a uniform distribution over the interval $[a, b]$. Using this notation, we describe the test systems used for simulations.

\begin{enumerate}
\item \emph{SoCal Distribution Circuit:} The sample industrial distribution system in Southern California has been previously reported in \cite{Masoud2011}. It has a peak load of approximately 11.3 MW and installed PV generation capacity of 6.4MW. We modified this circuit by removing the 30MW load at the substation bus (that represented other distribution circuits fed by the same substation) and simulated it with the parameters provided in Table \ref{table:opfParameters}. To scale the problem correctly, the problem was cast in per unit (p.u.) quantities using base values given in Table \ref{table:opfParameters}.

\item \emph{Random Test Circuits:} These circuits are generated using parameters typical of sparsely loaded rural circuits, as detailed in \cite{Schneider2008} and employed (with suitable modifications) in \cite{turitsyn11, YehGayme2012}. Around 15-60\% of the nodes are assumed to have 2 kW of PV capacity. The remaining parameters of these systems are described in Table  \ref{table:opfParameters}.
\end{enumerate}
The tests are run with both voltage and power-loss minimization as objective functions. The optimization results are summarized in Table \ref{table:opfResults}. For power-loss minimization, we always obtain a rank 1 optimal $W_*$. 

For voltage minimization, however, we obtain optimal solutions that violate the rank condition. In these cases, the heuristic approach is used to find a feasible point of $OPF$. We construct the solution based on the complex voltage $V_k = | V_k | e^{\ii \theta_k}$ at bus $k \in [n]$. For the heuristic algorithm, define

\bqn
x :=  ( |V_2|, |V_3|, \ldots, |V_n|, \theta_2, \theta_3, \ldots, \theta_n),
\eqn
and set the parameter $\gamma = +\infty$. In the examples studied, this approach always yields a feasible point within 5 iterations. From Table \ref{table:opfResults}, the values obtained for $\eta$ indicate that our algorithm finds a feasible point of $OPF$ with an objective value close to the theoretical optimum and hence performs well. A general bound on the performance of this heuristic technique, however, remains an open question.

\begin{table}
\begin{center}
\begin{tabular}{ | l || l | l |}
\hline
Test system
	& SoCal distribution circuit
	& Random radial networks\\
 \hline
 Number of nodes (n) & 47 & 50-150  \\
 \hline
Line impedances $(y_{ij})^{-1}$
	& \cite[Table 1]{Masoud2011}
	& $(0.33 + \ii 0.38) \Omega/km$, length = $(0.2km, 0.3km)$ \\
\hline
Voltage limits $\sqrt{\overline{W}_k}, \sqrt{\underline{W}_k}$
	& $1\pm 0.05$ p.u. at all nodes.
	& $1\pm 0.05$ p.u. at all nodes.\\
\hline
Real power demand $P_k^D$& \cite[Table 1]{Masoud2011} & $(0, 4.5 kW)$ \\
\hline
Reactive power demand $Q_k^D$ & Computed with $p.f. = (0.80, 0.98)$ lagging & $(0.2 P_k^D, 0.3 P_k^D)$ \\
\hline
Real power gen. limits
	& PV nodes: $\overline{P}_k^G = (0.2,1.0)$ times capacity,
	& PV nodes: $\overline{P}_k^G = (0, 2kW)$, \\
$\overline{P}_k^G$, $\underline{P}_k^G$
	& Substation node: $\overline{P}_1^G = 10MW$.
	& Substation node: $\overline{P}_1^G$ scaled with $n$.\\
	
	& At all nodes, $\underline{P}_k^G = 0$.
	& At all nodes, $\underline{P}_k^G = 0$. \\
 \hline
Reactive power gen. limits
	& $\overline{Q}_k^G= 0.3 \overline{P}_k^G, \ \ \underline{Q}_k^G = - 0.3 \overline{P}_k^G$ at all nodes.
	& $\overline{Q}_k^G= 0.3 \overline{P}_k^G, \ \ \underline{Q}_k^G = - 0.3 \overline{P}_k^G$ at all nodes.\\
 \hline
Base quantities & $P_{base} = 1MW$, $V_{base} = 12.35 kV (L-L).$ & $P_{base} = 1MW$, $V_{base} = 12.47 kV (L-L).$\\
\hline
\end{tabular}
\end{center}
 \caption{Circuit Parameters for SDR of $OPF$}
\label{table:opfParameters}
\end{table}

\begin{table}
\label{table:opfResults}
\begin{center}
\begin{tabular}{| l |c| c |c| c|}
\hline
Test system
	& \multicolumn{2}{|c|}{SoCal distribution circuit}
	& \multicolumn{2}{|c|}{Random radial networks}\\
\hline
 Minimize & Power-loss  & Voltage  & Power-loss  & Voltage \\
\hline
rank $W_*$ & 1 & $\geq 1$ &  1   &  $\geq1 $ \\
\hline
 Mean $\eta$ &  N/A & $1.8\%$   &   N/A & $0.5\%$\\
 \hline
Maximum $\eta$ &   N/A & $4.1\%$  & N/A &    $1.5\%$\\
 \hline
 \end{tabular}
\end{center}
\caption{Summary of simulation results}
\label{table:opfResults}
\end{table}



\section{Conclusion}
\label{sec:conc}
QCQP problems are generally non-convex and NP-hard. This paper proves that a certain class of 
 QCQP problems are solvable in polynomial-time. We have applied this result to the optimal power 
 flow problem and derived a set of conditions under which this nonconvex problem admits an 
 efficient solution. For problems that do not satisfy our sufficient conditions, we provide a 
 heuristic technique to find a feasible solution.   Simulations suggest  
 that this method often finds a near-optimal solution for the OPF problem.

\section{Appendix}
\label{app:proofs}

\subsection{Proof of lemma \ref{lemma:traceBound}:}
The result is restated here for convenience: If $H_1 \succeq 0$ and $H_2 \succeq 0$ are two $n \times n$ matrices,
\bqn
\tr (H_1  H_2) \geq  \rho_{\min}[H_1] \  \rho_{\max}[H_2].
\eqn
\begin{IEEEproof}
	Let the spectral decomposition of $H_2$ be
	\begin{align*}
	H_2 = \sum_{i=1}^n \rho_{i}  u_{i} u_{i}^{\mathcal{H}},
	\end{align*}
	where $\| u_i \| = 1$ for all $1 \leq i \leq n$. Then:
	\begin{align*}
	\tr(H_1 H_2) 
	& = \sum_{i =1}^n  \rho_{i}  \ (u_i^\mathcal{H} H_1 u_i) \\
	& \geq \sum_{i =1}^n  \rho_{i}  \ \rho_{\min}[H_1] \\
	& \geq \rho_{\min}[H_1] \ \rho_{\max}[H_2].
	\end{align*}
\end{IEEEproof}
\subsection{Proof of Theorem \ref{thm:mainResult} for $C \succeq 0$:}
The sketch of the proof is as follows. Perturb $RP$ so that the matrix in the objective function is positive definite. From our previous analysis this perturbed problem has a finite optimizer that has rank at most 1. Also, this perturbed problem can be solved in polynomial time. Using the solutions from the perturbed problems, we construct an optimum solution of $P$ in polynomial time.

In particular, consider the perturbed problems for $\delta > 0$:\\
\textbf{Perturbed primal problem $P^{(\delta)}$:}
\bqn
\label{QCQP:P}
\underset{x \in \mathbb{C}^{n}}{\text{minimize}}  & & x^{\mathcal{H}} (C + \delta \mathcal{I})  x
\\
\text{subject to:}
& &  x^{\mathcal{H}} C_k x  \leq b_k, \quad k \in \mathcal{K}. \label{eq:P.ineqC}
\eqn
\noindent
\textbf{Perturbed relaxed problem ${RP}^{(\delta)}$:}
\bqn
\underset{W \succeq 0}{\text{minimize}}  & & \tr [ (C + \delta \mathcal{I}) W ] \\
\text{subject to:}
& & \tr \left( C_k W \right) \leq {b}_k, \quad k \in \mathcal{K},
\eqn
where $\mathcal{I}$ is the $n \times n$ identity matrix. For any variable $z$ in $P/ RP$, let ${z}^{(\delta)}$ denote the corresponding variable in ${P}^{(\delta)}/ {RP}^{(\delta)}$. 

The matrix $C + \delta \mathcal{I}$ is positive definite for all $\delta > 0$. There exists ${W}_*^{(\delta)} \succeq 0$ that solves ${RP}^{(\delta)}$ and $\text{rank }{W}_*^{(\delta)} \leq 1$. Let the spectral decomposition of $W_*^{(\delta)}$ be $W_*^{(\delta)} = (x_*^{(\delta)})( x_*^{(\delta)})^{\mathcal{H}}$. From Lemma \ref{lemma:relaxedProblem}, we have
\bqn
p_*^{(\delta)} =  (x_*^{(\delta)})^{\mathcal{H}}(C + \delta \mathcal{I}) (x_*^{(\delta)}) = r_*^{(\delta)} = \tr [(C + \delta \mathcal{I})W_*^{(\delta)}] .
\eqn

The feasible regions of $P$ and $P^{(\delta)}$ are the same and are bounded. Then $x_*^{(\delta)}$ lies in a compact space, independent of $\delta$. Taking a sequence $\delta \to 0$, the corresponding sequence of $x_*^{(\delta)}$ has a convergent subsequence. Let the limit point of this subsequence be $\hat{x}$. Then $\hat{x}$ is feasible for $P$. We now show that $\hat{x}$ solves $P$ optimally.

The objective function of $P^{(\delta)}$ increases with $\delta$. Therefore,
\begin{align}
\label{eq:psd.1}
p_* \leq \hat{x}^\mathcal{H} C \hat{x} \leq p_*^{(\delta)} = r_*^{(\delta)}.
\end{align}
We wish to show that the first inequality is an equality. Suppose on the contrary $p_* < \hat{x}^\mathcal{H} C \hat{x}$. Let $x'_*$ be any finite optimizer of $P$. Then, 
\begin{align}
\label{eq:psd.2}
p_* =( x'_*)^\mathcal{H} C (x'_*) < \hat{x}^\mathcal{H} C \hat{x},
\end{align}
and we can choose a sufficiently small $\delta >0$, such that
\bqn
p_* + \delta  (x'_*)^\mathcal{H} ( x'_*) < \hat{x}^\mathcal{H} C \hat{x} \leq {r}_*^{(\delta)}.
\eqn
This follows from equations \eqref{eq:psd.1} and \eqref{eq:psd.2}. For this $\delta$, $[(x'_*) ( x'_*)^\mathcal{H} + (x_*^{(\delta)})( x_*^{(\delta)})^{\mathcal{H}}]/2$ is a feasible point of $RP^{(\delta)}$ and satisfies
\begin{align*}
	{r}_*^{(\delta)} 
	& \leq \tr \left[ (C + \delta \mathcal{I}) \left(\frac{(x'_*) ( x'_*)^\mathcal{H} + (x_*^{(\delta)})( x_*^{(\delta)})^{\mathcal{H}}}{2}\right) \right] \\
	& =  \frac{1}{2}{r}_*^{(\delta)}  + \frac{1}{2} \left[ p_* +  {\delta} (x'_*)^\mathcal{H} ( x'_*)  \right]\\
	& < {r}_*^{(\delta)}.
\end{align*}
This is a contradiction and hence $p_* = \hat{x}^\mathcal{H} C \hat{x}$. 

Now, we show how to use this perturbation technique to solve $P$ in polynomial time. Solve $RP$ to get $p_* = r_*$. If the optimizer $W_*$ of $RP$ has rank at most 1, compute $x_*$ from $W_*$ as in lemma 
\ref{lemma:relaxedProblem} then we have solved $P$ in polynomial time. If it does not satisfy the rank condition, then choose an arbitrary $\delta_0> 0 $ and solve $RP^{(\delta_0)}$
 in polynomial time to get the minimizer $W_*^{(\delta_0)}$ and the minimum $r_*^{(\delta_0)} = p_*^{(\delta_0)}$.
For any $\delta$ in $(0, \delta_0)$,
\begin{align}
\label{eq:convP.1}
p_* = r_* \leq \tr(C W_*^{(\delta)}) \leq p_*^{(\delta)}.
\end{align}
Also, $p_*^{(\delta)}$ is convex in $\delta$ \cite{yildirim2001sensitivity} and hence
\begin{align}
\label{eq:convP.2}
p_*^{(\delta)} \leq p_* + \frac{\delta}{\delta_0} \left(p_*^{(\delta_0)} - p_*\right).
\end{align}
Given $\zeta>0$, choose $\delta$ sufficiently small so that $ \frac{\delta}{\delta_0} \left(p_*^{(\delta_0)} - p_*\right) \leq \zeta$. For 
this $\delta$, solve $RP^{(\delta)}$ arbitrarily closely in polynomial time to get $W_*^{(\delta)}$ that has rank at most 1 and compute $x_*^{(\delta)}$. 
From equations \eqref{eq:convP.1} and \eqref{eq:convP.2}, $x_*^{(\delta)}$ satisfies
\bqn
\left | (x_*^{(\delta)})^\mathcal{H} C (x_*^{(\delta)}) - p_* \right | \leq \zeta.
\eqn
Thus $x_*^{(\delta)}$ is a feasible point of $P$ that has a value of the objective 
function within $\zeta$ of the theoretical optimum and we have shown a polynomial 
time algorithm to compute it.

This completes the proof of Theorem \ref{thm:mainResult}.


\subsection{Matrices involved in $OPF$:}

Here we compute the $(i, j)$-th entries of $\left\{ C_k, k \in \mathcal{K} \right\}$ for $OPF$. From \eqref{eq:defPhiPsi}, \eqref{eq:defPij}, \eqref{eq:defTij}, we have the following relations for $k \in [n]$, and
 $(p, q)$ and $(i,j)$ in $\mathcal{T}$:
\bq
[\Phi_k]_{ij} & = &  \begin{cases}
		\frac{1}{2} Y_{ij}  \, =\, \frac{1}{2} (-g_{ij} + \ii b_{ij}) &  \text{ if }k = i\\
		\frac{1}{2} {Y}_{ij}^{\mathcal{H}}  \,=\, \frac{1}{2} (-g_{ij} - \ii b_{ij}) &  \text{ if }k = j \\
		0 & \text{ otherwise},
		\end{cases}
\label{eq:Phi}
\eq
\bq
[\Psi_k]_{ij} & = &  \begin{cases}
		\frac{-1}{2\ii} Y_{ij}  \,=\, \frac{1}{2}(-b_{ij} - \ii g_{ij})  &  \text{ if }k = i \\
		\frac{1}{2\ii} {Y}_{ij}^{\mathcal{H}} \,=\, \frac{1}{2} (-b_{ij} + \ii g_{ij})  &  \text{ if }k = j \\
		0 & \text{ otherwise},
		\end{cases}
\label{eq:Psi}
\eq
\bq
[M^{pq}]_{ij} & = &  \begin{cases}
		g_{pq}  &  \text{ if } i=j=p \\
		\frac{1}{2} (-g_{pq} + \ii b_{pq})  &  \text{ if } (i, j) = (p, q) \\
		\frac{1}{2} (-g_{pq} - \ii b_{pq})  &  \text{ if } (i, j) = (q,p) \\
		0 & \text{ otherwise},
		\end{cases}
\label{eq:M}
\eq
\bq
[T^{pq}]_{ij} & = & \begin{cases}
		g_{pq}  &  \text{ if } i=j=p \text{ or } i=j=q \\
		-g_{pq}  &  \text{ if } (i, j) = (p,q) \text{ or } (i,j) = (q,p) \\
		0 & \text{ otherwise}.
		\end{cases}
\label{eq:T}
\eq

\bibliography{PowerbibNew}
\end{document}